\documentclass[11pt]{amsart}
\newtheorem{thm}[equation]{Theorem}
\newtheorem{pro}[equation]{Proposition}
\newtheorem{cor}[equation]{Corollary}
\newtheorem{lem}[equation]{Lemma}
\newtheorem{exa}[equation]{Example}

\newtheorem{DEF}[equation]{Definition}
\newtheorem{rem}[equation]{Remark}

\def\iso{\hbox{iso}}

\def\andd{\quad\hbox{and}\quad}
\def\ind{\hbox{ind}}

\def\v{{\mathcal V}}
\def\u{{\mathcal U}}
\def\vd{\dot{\mathcal V}}

\def\fm{(\cdot,\cdot)}
\def\a{\alpha}

\def\w{{\mathcal W}}

\def\wdot{\dot{\mathcal W}}
\def\sub{\subseteq}
\def\rd{\dot{R}}

\def\rt{\tilde{R}}

\def\lam{\lambda}
\def\Lam{\Lambda}

\def\1k{\frac{1}{k}}

\def\la{\langle}
\def\ra{\rangle}

\def\d{\delta}

\def\b{\beta}

\def\qed{\hfill$\Box$}

\def\sg{\sigma}

\def\rtimes{R^{\times}}
\def\gg{{\mathcal G}}

\def\hh{{\mathcal H}}

\def\gs{\gg^{\sg}}
\def\hs{\hh^{\sg}}

\def\ta{t_\a}

\def\ad{\hbox{ad}}

\def\gt{\tilde{\gg}}

\def\gc{\gg_{c}}

\def\riso{R_{\hbox{iso}}}

\def\rniso{R_{\hbox{niso}}}

\def\i{{\mathcal I}}

\def\bbbc{{\mathbb C}}
\def\bbbz{{\mathbb Z}}
\def\bbbr{{\mathbb R}}

\def\hd{\dot{\hh}}

\def\dd{\mathcal D}

\def\aff{\hbox{Aff}}
\def\hht{\tilde{\mathcal H}}

\def\gfh{(\gg,\fm,\hh)}
\def\ga{\gg_{\a}}
\def\hdi{\dot{\hh}_{i}}

\def\span{\hbox{span}}
\def\vdi{\dot{\v}_i}
\def\rdi{\dot{R}_i}
\def\hi{\hh_i}
\def\da{\dot{\a}}
\def\hd{\dot{\hh}}
\def\wdot{\dot{\w}}

\def\proof{\noindent{\bf Proof. }}
\def\vb{\bar{\v}}
\def\ab{\bar{\a}}
\def\bb{\bar{\b}}
\def\rb{\bar{R}}
\def\lrd{\la\rd\ra}
\def\db{\dot{\b}}
\def\ll{\mathcal L}
\def\zgc{{\mathcal Z}(\gc)}
\begin{document}
\markboth{GENERALIZED REDUCTIVE LIE ALGEBRAS}{S. AZAM}

\centerline{\bf GENERALIZED REDUCTIVE LIE ALGEBRAS,}
\centerline{connections with} \centerline{extended affine Lie
algebras and Lie tori} \vspace{.5cm}\centerline{SAEID
AZAM\footnote[1]{Research number 810821, University of Isfahan,
Iran. This research was in part supported by a grant from IPM (No.
82170011).}}

\vspace{.5cm}\centerline{Department of Mathematics, University of
Isfahan,} \centerline{P.O.BOX 81745, Isfahan, Iran,\&}
\centerline{ Institute for Theoretical Physics and Mathematics
(IPM)}

%\footnote{1991 Mathematics Subject classification. Primary 17B65,
%Secondary 17B67, 17B40.}
%\footnotetext[1]{Research number
%810821, University of Isfahan, Iran.} \footnote[2]{This research
%was in part supported by a grant from IPM.}
%\date{}
%\maketitle
%\begin{abstract}

\vspace{1cm}\centerline{\bf ABSTRACT}
 We investigate a class of Lie algebras which we
 call {\it generalized reductive Lie
 algebras}. These are generalizations of  semi-simple,
reductive, and affine Kac-Moody Lie algebras. A generalized
reductive Lie algebra which has an irreducible root system is said
to be {\it irreducible} and we note that this class of algebras
have been under intensive investigation in recent years. They have
also been called {\it extended affine Lie algebras}. The larger
class of generalized reductive Lie algebras has not been so
intensively investigated. We study them in this paper and note
that one way they arise is as fixed point subalgebras of finite
order automorphisms. We show that the core modulo the center of a
generalized reductive Lie algebra is a direct sum of centerless
Lie tori. Therefore one can use the results known about the
classification of centerless Lie tori to classify the cores modulo
centers of generalized reductive Lie algebras.

\vspace{5mm} \setcounter{section}{-1}
\section{\bf Introduction}
In 1990 H\o egh-Krohn and B Torresani [HK-T] introduced a new
interesting class of Lie algebras over field of complex numbers,
called {\it quasi simple Lie algebras} by proposing a system of
fairly natural and not very restrictive axioms. These Lie algebras
are characterized by the existence of a symmetric nondegenerate
invariant bilinear form, a finite dimensional Cartan subalgebra, a
discrete root system which contains some nonisotropic roots, and
the ad-nilpotency of the root spaces attached to non-isotropic
roots. As it will appear from the sequel, these algebras are
natural generalizations of reductive Lie algebras, and affine
Kac--Moody Lie algebras. For this reason and other reasons
indicated in the introduction of the paper [AABGP] we call this
class of Lie algebras {\it generalized reductive Lie algebras}
(GRLA for short).

In [HK-T], the authors extract some basic properties of GRLAs from
the axioms, but for the further study of such Lie algebras they
assume the irreducibility of the corresponding root systems.
Namely, a GRLA is called {\it an extended affine Lie algebra}
(EALA for short) if the set of non-isotropic roots is
indecomposable and isotropic roots are non-isolated (see
Definition \ref{qssla} for terminology). We note that EALAs have
been under intensive investigation in recent years, however the
more general class of generalized reductive Lie algebras has not
been so intensively investigated.

In [AABGP] the axioms for an EALA are introduced in steps in such
a way that the power of each axiom is clearly shown before
introducing the next one. This in particular provides a good
framework for the study of Lie algebras which satisfy only a part
of the axioms. In Section 1, we follow the same steps as in
[AABGP, Chapter I] to obtain the basic structural properties of a
GRLA. However, it turns out that some of the results in [AABGP,
Chapter I] can be proved using fewer axioms than those in
[AABGP](see Proposition \ref{w6}(i)).

In Section 2, we get more information about the structure of a
GRLA $\gg$ by decomposing its root system $R$ into a finite union
of indecomposable subroot systems, and then corresponding to each
subroot system we construct an indecomposable generalized
reductive subalgebra of $\gg$. More precisely, we show that up to
some isolated spaces, a GRLA is a finite sum of certain
indecomposable generalized reductive subalgebras and an abelian
subalgebra, with trivial Lie bracket between distinct summands (on
the level of core). In particular, if there are no isolated root
spaces (that is if $\gg$ satisfies part (b) of axiom GR6 of a
GRLA), the structure of $\gg$ can be thought of a generalization
of a reductive Lie algebra. In fact when the nullity is zero,
$\gg$ is nothing but a reductive Lie algebra (see Corollary
\ref{main2} for details). The main result of this section is that
the core modulo center of a GRLA $\gg$ is isomorphic to a direct
sum of the cores modulo centers of some indecomposable generalized
reductive ideals of $\gg$. When the nullity is less than or equal
two, this result can be read as the core of $\gg$ modulo its
center is a direct some of the cores modulo centers of some
extended affine Lie algebras (see Theorem \ref{main}(iv) and
Corollary \ref{main2}(i)).

In Section 3, the main section, we show that the core modulo
center of an indecomposable GRLA is a centerless Lie torus, and
therefore the core modulo center of a GRLA is a direct sum of
centerless Lie tori (Theorem \ref{classification} and Corollary
\ref{final}). Therefore, one can use the results of [BGK], [BGKN],
[AG], [Y], [ABG] and [AFY] regarding the classification of
centerless Lie tori to classify the cores modulo centers of GRLAs
for types which the classification is achieved. In principle, the
classification of centerless Lie tori is done for all types except
type $BC_2$. See also [A2, Proposition 1.28], [AG, Proposition
1.28] and [N2, Theorem 6] for the relation between an EALA and its
core modulo center. For a deep study of EALAs and their root
systems we refer the reader to [S], [BGK], [BGKN], [AABGP], [AG],
[ABG], [A1], [A3]. Also see [N1-2] and [AKY] for some new classes
of Lie algebras which are closely related to EALAs.

In Section 4, we give several examples of GRLAs and we show some
methods of constructing new GRLAs from old ones. In particular, it
is shown that GRLAs arise as the fixed point subalgebras of finite
order automorphisms.

The author would like to thank the referee for pointing out a gap
in an early version of this work, as well as for many helpful
suggestions. He also would like to thank Professors, B. Allison,
S. Berman, K. H. Neeb and A. Pianzola for some helpful
discussions. A discussion with Professor K. H. Neeb was led to
Example \ref{exam3}. A part of this work was completed while the
author was visiting the Department of Mathematical Sciences at the
University of Alberta, thanks for their hospitality.

\vspace{5mm}
\section{\bf Generalized reductive Lie algebras}
Let $\gg$ be a Lie algebra over the field of complex numbers, let
$\hh$ be a subalgebra of $\gg$ and
$\fm:\gg\times\gg\rightarrow{\bbbc}$ be a bilinear form on $\gg$.
Consider the following axioms for the triple $\gfh$:

GR1. The form $\fm$ is symmetric, nondegenerate and invariant on
$\gg$.

GR2. $\hh$ is a nontrivial finite dimensional abelian subalgebra
which is self-centralizing and $\ad(h)$ is diagonalizable for all
$h\in\hh$.

According to GR2 we have a vector space decomposition
$\gg=\oplus_{\a\in \hh^\star}\gg_{\a}$, where
$$\ga=\{x\in\gg\mid[h,x]=\a(h)x\hbox{ for all }h\in\hh\}.$$
The set $R=\{\a\in\hh^\star\mid\gg_\a\not=\{0\}\}$ is called the
{\it root system} of $\gg$.

From GR1-GR2 it follows that
\begin{equation}\label{w1}
\gg_0=\hh,\quad 0\in R,
\end{equation}
and
\begin{equation}\label{w2}
(\ga,\gg_\b)=\{0\}\hbox{ unless }\a+\b=0.
\end{equation}
In particular,
\begin{equation}\label{w3}
R=-R,
\end{equation}
and the form restricted to $\hh$ is nondegenerate. For
$\a\in\hh^\star$ let $\ta$ be the unique element in $\hh$ which
represents $\a$ via the form. Then for any $\a\in R$,
\begin{equation}\label{w4}
[\ga,\gg_{-\a}]={\bbbc}t_\a.
\end{equation}
Transfer the form to $\hh^\star$ through
\begin{equation}\label{ww4}
(\a,\b):=(t_\a,t_\b)\hbox{ for }\a,\b\in\hh^\star.
\end{equation}
Let
$$
R^{\times}=\{\a\in R\mid (\a,\a)\not=0\}\andd R^0=\{\a\in R\mid
(\a,\a)=0\}.
$$
Elements of $\rtimes$ $(R^0$) are called non-isotropic (isotropic)
roots of $R$. The next axioms are as follows:

GR3. For any $\a\in R^\times$ and $x\in\ga$, $\ad_{\gg}(x)$ acts
locally nilpotently on $\gg$.

GR4. $R$ is a discrete subset of $\hh^\star$.

GR5. $\rtimes\not=\emptyset$.

\begin{DEF}\label{qssla}
{\rm A triple $\gfh$ satisfying axioms GR1-GR5 is called a {\it
generalized reductive Lie algebra} (GRLA for short). We call a
generalized reductive Lie algebra {\it indecomposable} if it
satisfies

GR6a. $\rtimes$ is indecomposable, that is $\rtimes$ is not a
disjoint union of two of its nonempty subsets which are orthogonal
with respect to the form.

We call a GRLA {\it non-singular} if it satisfies

GR6b. For $\sg\in R^0$, there  exists $\a\in\rtimes$ such that
$\a+\sg\in R$, that is isotropic roots of $R$ are non-isolated.

Finally, a GRLA is called an {\it extended affine Lie algebra}
(EALA for short) if it satisfies GR1-GR6. When there is no
confusion, we simply write $\gg$ instead of $\gfh$. The {\it core}
of a GRLA $\gg$ is by definition, the subalgebra $\gc$ of $\gg$
generated by root spaces corresponding to non-isotropic roots. It
follows that $\gc$ is a perfect ideal of $\gg$. A GRLA is called
tame if $\gc$ contains its centralizer in $\gg$.}
\end{DEF}

\begin{rem}\label{w5}
{\rm (i) It follows from axioms GR1-GR2 and GR6b that
$\rtimes\not=\emptyset$, so the axiom GR5 is redundant for a
non-singular GRLA.

(ii) It is easy to see that a GRLA $\gg$ is tame if and only if
$C_{\gg}(\gc)=\{x\in\gg\mid (x,\gc)=\{0\}\}$. The proof of [ABP,
Lemma 3.62] shows that a tame GRLA is non-singular.

(iii) Semisimple Lie algebras, finite dimensional reductive Lie
algebras and a direct sum of EALAs are examples of non-singular
GRLAs.  Heisenberg Lie algebras (with derivation added) satisfy
axioms GR1-GR4, however they are not GRLAs as $\rtimes=\emptyset$.
It is shown in [ABY, Section 3] that the fixed point subalgebra of
an EALA under a finite order automorphism satisfies GR1-GR4.}
\end{rem}

{\it From now on we assume that $\gfh$ is a GRLA with the
corresponding root system $R$}. So we may use all the results in
[AABGP] which are obtained by axioms GR1-GR5. Let us state from
[AABGP] here some of the important properties of $\gg$ which will
be of use in the sequel. We emphasize on the particular axioms
which are use in the proof of each result.

It is shown in [AABGP, I.(1.18)] that if $\gg$ satisfies GR1-GR2,
then for $\a\in\rtimes$ there exist $e_a\in\gg_\a$ and
$f_\a\in\gg_{-\a}$ such that
\begin{equation}\label{sl2}
(e_\a,h_\a:=[e_\a,f_\a],f_\a)\hbox{ is a
}\mathfrak{sl}_2\hbox{-triple,} \end{equation}
 that is the
$\bbbc$-span of $\{e_a,h_\a,f_{\a}\}$ is a Lie subalgebra of $\gg$
isomorphic to $\mathfrak{sl}_2(\bbbc)$. Note that
$$
h_\a=\frac{2t_\a}{(\a,\a)}.
$$

If $\u$ is a vector space equipped with a bilinear form, let us
write
$$
\a^\vee:=\frac{2\a}{(\a,\a)}\hbox{ for }\a\in\u\hbox{ with
}(\a,\a)\not=0.
$$

\begin{thm}\label{aabgp}[AABGP, Theorem I.1.29]
Let $\gg$ satisfy GR1-GR3 and $\a\in\rtimes$. Then

(a) For $\b\in R$, we have $(\b,\a^\vee)\in{\bbbz}$.

(b) For $\b\in R$, $\b-(\b,\a^\vee)\a\in R$.

(c) If $k\in {\bbbc}$ and $k\a\in R$, then $k=0,\pm 1$.

(d) $\dim \ga=1$.

(e) For any $\b\in R$, there exist two non-negative integers $u,d$
such that for any $n\in{\bbbz}$ we have $\b+n\a\in R$ if and only
if $-d\leq n\leq u$. Moreover, $d-u=(\b,\a^\vee)$.
\end{thm}

The statement of part (i) of the following proposition is the same
as [AABGP, Proposition I.2.1], however the proof given here is
different, as we are not allowed to use axiom GR6. We even do not
use GR4 in the proof.

\begin{pro}\label{w6} (i) Let $\gg$ satisfy GR1-GR3. Then $(R,R^0)=\{0\}$.

(ii) $\d\in R^0$ is isolated if and only if $\gg_\d\sub
C_{\gg}(\gc)$.
\end{pro}

\proof (i) First let $\a\in\rtimes$ and $\d\in R^0$. Suppose to
the contrary that $(\a,\d)\not=0$. By [AABGP, Lemma I.1.30],
$\a+n\d\in R$ for sufficiently large $n$, and it is clear that
$\a+n\d\in\rtimes$ except at most for one $n$. But then for
suitable $n$ we have
$$
\frac{2(\a+n\d,\d)}{(\a+n\d,\a+n\d)}=
\frac{2(\a,\d)}{(\a,\a)+2n(\a,\d)}\notin{\bbbz}$$ which
contradicts part (a) of Theorem \ref{aabgp}.

Next let $\d,\eta\in R^0$. We must show $(\d,\eta)=0$. If not,
then $\eta+\d$ and $\eta-\d$ are non-isotropic and are not
orthogonal to $\d$, $\eta$. Therefore by the previous step,
$\eta\pm\d\not\in R$. So we get a contradiction if we show that
$$(\eta,\d)\not=0\Rightarrow \eta+\d\in R\hbox{ or
}\eta-\d\in R.
$$ Suppose $(\eta,\d)\not=0$ and $\eta-\d\not\in
R$. Choose $x_\d\in\gg_{\d}$ and $x_{-\d}\in\gg_{-\d}$ such that
$[x_{-\d},x_{\d}]=t_{\d}$. Take any $0\not=x_{\eta}\in\gg_{\eta}$.
Then using the Jacobi identity, we have
$$[x_{-\d},[x_{\d},x_{\eta}]]=(\eta,\d)x_{\eta}\not=0.
$$
Thus $[x_{\d},x_\eta]\not=0$, and so $\eta+\d\in R$.

(ii) Suppose first that $\d\in R^0$ is not isolated, that is there
exists $\a\in R^\times$ such that $\a+\d\in R$. We must show that
$[\gg_\d,\gc]\not=\{0\}$. Since $\gg_\a\sub\gc$, it is enough to
show that $[\gg_\d,\ga]\not=\{0\}$. Consider the non-negative
integers $u,d$ appearing in the $\a$-string through $\d$, as in
part (e) of Theorem \ref{aabgp}. We have
 $d=u$ as by
part (i), $(\d,\a)=0$. Let $0\not=x\in\gg_{\d-d\a}$. Then
$[x,\gg_{-\a}]\sub\gg_{\d-(d+1)\a}=\{0\}$. By [AABGP, Lemma
I.1.21] and part (d) of Theorem \ref{aabgp}, for any
$0\not=z\in\gg_\a$, we have
$$
(\ad z)^N(x)\not=0\quad\hbox{but}\quad(\ad z)^{N+1}(x)=0,
$$
where $N=2(\d-d\a,-\a)/(\a,\a)=2d$. Since $\a+\d\in R$, we have
$d=u\geq 1$, so $d+1\leq 2d=N$. Thus $[z,(\ad z)^d(x)]=(\ad
z)^{d+1}(x)\not=0$. But $(\ad z)^d(x)\in\gg_\d$ and so
$[\gg_\d,\ga]\not=\{0\}$. Conversely, if $\d$ is isolated then
$\a+\d\not\in R$ for all $\a\in\rtimes$. Thus
$[\gg_\d,\ga]=\{0\}$.\hfill\qed\vspace{3mm}

\begin{rem} {\rm (i) According to part (ii) of Proposition \ref{w6},
axiom GR6b is equivalent to

GR6b$'$: For any $\d\in R^0$, the root space $\gg_\d$ is not
contained in the centralizer of the core.

(ii) The proof of part (ii) of Proposition \ref{w6} and [AG, Lemma
1.3] show that for any $\a\in\rtimes$ and $\b\in R$, $\a+\b$ is a
root if and only if $[\ga,\gg_\b]\not=\{0\}$. In particular, if
$\a+\b\in\rtimes$, then $[\ga,\gg_\b]=\gg_{\a+\b}$. The proof of
Proposition \ref{w6}(ii) is in fact a modified version of a
standard $\mathfrak{sl}_2$-argument.}
\end{rem}

Define an equivalence relation on $\rtimes$ by saying that two
roots $\a$ and $\b$ in $\rtimes$ are related if and only if there
is a sequence of roots $\a_0=\a,\a_1,\ldots,\a_t=\b$ in $\rtimes$
such that $(\a_i,\a_{i+1})\not=0$ for $0\leq i\leq t-1$. This
determines a partition $\rtimes=\cup_{i\in I}\rtimes_i$ where
$\rtimes_i$'s are indecomposable and $(\rtimes_i,\rtimes_j)=0$ for
$i\not=j$. In particular, if $\v$ is the real span of $R$ and
$\v_i$ is the real span of $\rtimes_i$ then
\begin{equation}\label{indecom}
(\v_i,\v_j)=0\quad\hbox{for}\quad i\not=j,
\end{equation}
and
\begin{equation}\label{indecom1}
\v=\sum_{i\in I}\v_i+\span_{\bbbr}R^0.
\end{equation}
%Note that by Proposition \ref{w6}(i)
%$$
%\span_{\bbbr}R^0\sub\v^0.
%$$

Fix $i\in I$. It follows from part (a) of Theorem \ref{aabgp} and
indecomposability of $\rtimes_i$ (see [AABGP, I.\S 2]) that there
exist nonzero scalars $c_i\in\bbbc$ such that the form on $\v_i$
defined by
\begin{equation}\label{scaled}
\fm_i:=c_i\fm,
\end{equation}
is real valued and
\begin{equation}\label{pos0}
(\gamma,\gamma)_i>0\quad\hbox{for some}\quad \gamma\in\rtimes_i.
\end{equation}
If (\ref{pos0}) does not hold for all $\gamma\in\rtimes_i$, then
it follows from indecomposability of $\rtimes_i$ that there are
roots $\a,\b\in\rtimes_i$ such that $(\a,\a)_i>0$, $(\b,\b)_i<0$
and $(\a,\b)_i\not=0$. But as the proof of [AABGP, Lemma I.2.3]
suggests this leads to the existence of a complex simple Lie
algebra of dimension $6$ which is absurd. Thus
\begin{equation}\label{pos1}
(\a,\a)_i>0\quad\hbox{for all}\quad \a\in\rtimes_i.
\end{equation}
Also note that for $\a,\b\in\v_i$ with $(\a,\a)_i\not=0$, we have
\begin{equation}\label{pos3}
(\b,\a^\vee)_i:=\frac{2(\b,\a)_i}{(\a,\a)_i}=(\b,\a^\vee).
\end{equation}

The proof of the following lemma is almost the same as [AABGP,
Lemma I.2.6], however we have to be careful for one part in which
the $\a$-string through $\b$ is used. We provide the proof here
with the necessary modifications.
\begin{lem}\label{pos5}
For $\a,\b\in\rtimes$, we have $-4\leq (\b,\a^\vee)\leq 4$.
\end{lem}

\proof Since $(\rtimes_i,\rtimes_j)=\{0\}$ for $i\not= j$ and
$(R,R^0)=\{0\}$ by Proposition \ref{w6}, we may assume that
$\a,\b\in\rtimes_i$ for some $i$. Since
$(\b,\a^\vee)=(\b,\a^\vee)_i$ we may replace $\fm$ with $\fm_i$.
So it is enough to show that if $(\a,\b)_i<0$ then
$(\b,\a^\vee)_i\geq -4.$ Suppose to the contrary that
$(\a,\b)_i<0$ but $a=(\b,\a^\vee)\leq -5$. Let $b=(\a,\b^\vee)$.
From (\ref{pos1}) we have $b\leq -1$. By Theorem \ref{aabgp}(e)
all elements of the string
$$
\b-d\a,\ldots,\b-\a,\b,\b+\a,\ldots,\b+u\a,
$$
are elements of $R$, where $d-u=a\leq -5$. Thus $u\geq 5$. In
particular $\b+2\a\in R$. If $(\b+2\a,\a)_i=0$ then
$(\b,\a)_i=-2(\a,\a)_i$. So
$$
-5\geq a=\frac{2(\b,\a)_i}{(\a,\a)_i}=
\frac{-4(\a,\a)_i}{(\a,\a)_i}=-4$$ which is absurd. Thus
$(\b+2\a,\a)_i\not=0$ and so $\b+2\a\in\rtimes_i$. Then
$$(\b+2\a,\b+2\a)_i=\frac{2(\b,\a)_i}{a}(a+4)<0.$$
This contradicts (\ref{pos1}).\qed

The restriction of the form $\fm$ to $\v\times\v$ defines a
symmetric bilinear map
$$\fm:\v\times\v\longrightarrow\bbbc.$$
Let
$$\v^0:=\{v\in\v\mid (v,\v)=\{0\}\},$$
be the radical of this map. Let $\vb=\v/\v^0$ and consider the
canonical map $\bar{\;}:\v\rightarrow\vb$. We have $\vb\not=\{0\}$
as we have assumed that $\rtimes\not=\emptyset$. Define on the
real vector space $\vb\times\vb$ a complex valued symmetric
bilinear map
\begin{equation}\label{barform}
\fm:\vb\times\vb\rightarrow\bbbc \end{equation} by
$$(\ab,\bb):=(\a,\b)\quad\hbox{for}\quad\a,\b\in\v.$$
Then $\fm$ is nondegenerate on $\vb$. Moreover, if $\vb_i$ is the
image of $\v_i$ under the map $\bar{\;}$, then the restriction of
$\fm$ to $\vb_i\times\vb_i$ is also nondegenerate and for
$\ab,\bb\in\vb_i\setminus\{0\}$,
\begin{equation}\label{z}
(\bb,\ab^\vee)_i:=\frac{2(\bb,\ab)_i}{(\ab,\ab)_i}=(\bb,\ab^\vee)
\in\bbbz.
\end{equation}

Set
$$
\bar{R}=\{\ab\mid\a\in
R\}\andd\rb_i=\{\ab\mid\a\in\rtimes_i\}\cup\{0\}.
$$

\begin{lem}\label{pos6}
%\begin{equation}\label{pos6}
$\rb$ is finite.
%\end{equation}
\end{lem}

\proof Consider a basis $\{\ab_1,\ldots,\ab_\ell\}\sub\rb$ of
$\vb$. If $\ab_{j}\in \rb_{i_j}$ for some $i_j\in I$, then by
(\ref{pos1}) the term $c_{i_j}(\ab_j,\ab_j)$ is a nonzero real
number and so
$$\{\frac{2\ab_1}{c_{i_1}(\ab_1,\ab_1)},\ldots,\frac{2\ab_\ell}{c_{i_\ell
}(\ab_\ell,\ab_\ell)}\}$$ is also a basis of $\vb$. Now for
$\bb\in\vb$ define
$$
\varphi(\bb)=\big(c_{i_1}(\bb,\frac{2\ab_1}{c_{i_1}(\ab_1,\ab_1)}),\ldots,
c_{i_\ell}(\bb,\frac{2\ab_{\ell}}{c_{i_\ell}(\ab_\ell,\ab_\ell)})\big).
$$
Note that for each $j$, if $(\bb,\ab_j)\not=0$, then
$\bb\in\rb_{i_j}$ and so by (\ref{z}),
$\varphi(\bb)\in\bbbz^\ell$. Now from Lemma \ref{pos5} we see that
$\varphi(\rb)$ has at most $9^\ell$ elements. Since the bilinear
map (\ref{barform}) on $\vb$ is nondegenerate, $\varphi$ is one to
one and so $\rb$ is finite.\qed

We have from Lemma \ref{pos6} that the index set I is finite, say
$I=\{1,\ldots,k\}$, and so
\begin{equation}\label{well}
\rtimes=\cup_{i=1}^{k}\rtimes_i\andd\v=\sum_{i=1}^{k}\v_i+\hbox{span}_{\bbbr}R^0,
\end{equation}
where $\span_{\bbbr}R^0\sub\v^0$.

For $\a\in\rtimes$ define $w_\a:\v\rightarrow\v$ by
$$
w_\a(\b)=\b-(\b,\a^\vee)\a,\quad(\b\in\v).
$$
By Theorem \ref{aabgp}(b), $w_\a(R)\sub R$. Since $w_\a$ preserves
the form $\fm$, we have from (\ref{indecom}) that for
$\a\in\rtimes_i$,
$$
w_\a(\rtimes_i)\sub\rtimes_i\andd w_\a(\v_i)\sub\v_i.
$$
In a similar manner, for $\ab\in\rb_i\setminus\{0\}$, we can
define $w_{\ab}$ on $\vb_i$. Then
\begin{equation}\label{bar}
w_{\ab}(\bb)=\overline{w_\a(\b)}\andd w_{\ab}(\rb_i)\sub\rb_i.
\end{equation}

\begin{lem}\label{pos7}
(i) The symmetric bilinear form $\fm_i$ on $\vb_i$ is positive
definite.

(ii) The symmetric bilinear form $\fm_i$ on $\v_i$ is positive
semidefinite.
\end{lem}

\proof Clearly (ii) is a immediate consequence of (i). By Lemma
\ref{pos6}, $\rb_i$ is finite and $\rb_i\setminus\{0\}$ is
indecomposable with respect to the form $\fm_i$, as $\rtimes_i$ is
indecomposable. Now follow the proof of [AABGP, Theorem I.2.14]
with $\rb_i$ in place of $\rb$ and $\vb_i$ in place of $\vb$.\qed

We now would like to put together the forms $\fm_i$ to obtain a
positive semidefinite bilinear form on $\v$. Using (\ref{well}) we
can write each element $\a\in\v$ in the form
\begin{equation}\label{pos8}
\a=\sum_{i=1}^{k}\a_i+\d_\a\quad\hbox{where}\quad\a_i\in\v_i,\;\;
\d_\a\in\span_{\bbbr}R^0\sub\v^0.
\end{equation}
\begin{lem}\label{form}
Let $\a=\sum_{i=1}^{k}\a_i+\d_\a$ and
$\b=\sum_{i=1}^{k}\b_i+\d_\b$ be two elements of $\v$ in the form
(\ref{pos8}) and define
\begin{equation}\label{newform}
(\a,\b)'=\sum_{i=1}^{k}(\a_i,\b_i)_i.
\end{equation}
Then $\fm'$ is a well-defined real valued positive semidefinite
symmetric bilinear form on $\v$. Moreover,
\begin{equation}\label{form1}
(\b,\a^\vee)':=\frac{2(\b,\a)'}{(\a,\a)'}=(\b,\a^\vee)\quad\hbox{for}\quad\a,\b\in\rtimes.
\end{equation}
\end{lem}

\proof To show that $\fm'$ is well-defined consider $\a$ and $\b$
as in the statement and let $\a=\sum_{i=1}^{k}\a'_i+\d'_{\a}$ and
$\b=\sum_{i=1}^{k}\b'_i+\d'_{\b}$ be other expressions of $\a$ and
$\b$ in the form (\ref{pos8}). We must show that
$$
\sum_{i=1}^{k}(\a_i,\b_i)_i=\sum_{i=1}^{k}(\a'_i,\b'_i)_i.
$$
Now using (\ref{indecom}) and the fact that
$\d_\a,\d_\b,\d'_\a,\d'_\b$ are isotropic, we have
\begin{eqnarray*}
(\a_i,\b_i)_i
&=&c_i(\a-\sum_{j\not= i}\a_j-\d_\a,\b_i)\\
&=&c_i(\sum_{j=1}^{k}\a'_j+\d_\a'-\sum_{j\not=i}\a_j-\d_\a,\b_i)\\
&=&c_i(\a'_i,\b_i)\\
&=&c_i(\a'_i,\b-\sum_{j\not=i}\b_j-\d_\b)\\
&=&c_i(\a'_i,\sum_{j=1}^{k}\b'_j+\d_\b'-\sum_{j\not=i}\b_j-\d_\b)\\
&=&c_i(\a'_i,\b'_i)=(\a'_i,\b'_i)_i.
\end{eqnarray*}
This proves that the form $\fm'$ is well-defined. Now since
$\fm_i$ is positive semidefinite on $\v_i$, for each $i$, it is
clear that $\fm'$ is positive semidefinite on $\v$ and that
(\ref{form1}) holds.\qed

We have from Lemmas \ref{pos7}, \ref{pos6}, (\ref{z}) and
(\ref{bar}) that
\begin{equation}\label{finite}
\rb_i\hbox{ is an irreducible finite root system in }\vb_i.
\end{equation}

Lemma \ref{form} together with other properties which we have seen
about $\v$ and $R$ lead us to state the following definition.

\begin{DEF}\label{EARS}
{\rm Let $\v$ be a nontrivial finite dimensional real vector space
with a nontrivial positive semidefinite symmetric bilinear form
$(.,.)$ and let $R$ be a subset of $\v$. We say $R$ is a {\it
generalized reductive root system} (GRRS for short) in $\v$ if $R$
satisfies the following 5 axioms:

%\noindent (R1) $0\in R$,

\noindent (R1) $R=-R$,

\noindent (R2) $\;R \hbox{ spans }\v$,

\noindent (R3) $\;$ $R$ is discrete in $\v$,

\noindent (R4) If $\a\in\rtimes$ and $\b\in R$, there exist two
non-negative integers $u,d$ such that for any $n\in{\bbbz}$ we
have $\b+n\a\in R$ if and only if $-d\leq n\leq u$. Moreover,
$d-u=(\b,\a^\vee)$.

\noindent (R5) $\;\a\in R^{\times} \Rightarrow 2\a\not\in R$.

We call the GRRS $R$ {\it non-singular} if it satisfies:

\noindent (R6) $\;$ for any $\d\in R^{0}$, there exists $\a\in
R^{\times}$ such that $\a+\d\in R$. We say a root satisfying this
condition is {\it nonisolated} and call isotropic roots which do
not satisfy this {\it isolated}.

The GRRS $R$ is called {\it indecomposable} if it satisfies:

\noindent (R7)$\;$ $R^{\times}$ cannot be decomposed into a
disjoint union of two nonempty subsets which are orthogonal with
respect to the form.

A non-singular indecomposable GRRS $R$ is known at the literature
as {\it an extended affine root system} (EARS for short). One may
also call it an {\it irreducible} GRRS. The {\it nullity} of a
GRRS $R$ is defined to be the dimension of the real span of
$R^0$.}
\end{DEF}

Since the form in the definition of a GRRS is nontrivial we have
from (R2) that $\rtimes\not=\emptyset$. Then it follows from this,
(R1) and (R4) that $0\in R$. The root system $R$ of a
(non-singular) GRLA $\gg$ is a (non-singular) GRRS. In fact the
existence of a nontrivial positive semi-definite bilinear form was
shown in Lemma \ref{form}, and by (\ref{w1}), (\ref{w3}), GR4 and
Theorem \ref{aabgp}(d)-(e) axioms (R1)-(R5) also hold. From Remark
\ref{w5}(ii) we know that the root system of a tame GRLA is a
non-singular GRRS.  We define the {\it nullity} of a GRLA to be
the nullity of its root system.

For a GRRS $R$ we set
$$
%\begin{array}{l}
\riso=\{\d\in R^0\mid\d+\a\not\in R\hbox{ for any
}\a\in\rtimes\},\hbox{ and } \rniso=R^0\setminus\riso.
$$
That is $\riso$ ($\rniso$) is the set of isolated (non-isolated)
roots of $R$. So $R$ is non-singular if and only if
$\riso=\emptyset$.

\vspace{5mm}
\section{\bf Intrinsic decomposition of a GRLA}
\setcounter{equation}{0}
 Let $\gfh$ be a GRLA with root system
$R$. Let $\v$, $\bar{\v}$, $\v_i$ and $\vb_i$ be as in Section 1.
Let $\fm$ be the form on $\hh^\star$, defined by (\ref{ww4}),
restricted to $\v$. Fix a real valued positive semidefinite
symmetric bilinear form $\fm_i$ on $\v_i$, as in (\ref{scaled}).
Let $\fm$ be the real positive semidefinite symmetric bilinear
form on $\v$ defined by (\ref{newform}). Then
$$\fm'_{\mid_{\v_i}}=\fm_i=c_i\fm.$$
We also have from Lemma \ref{pos7}(i) that the form $\fm_i$ on
$\vb_i$ is real-valued and positive definite, and $\rb_i$ is an
irreducible finite root system in $\vb_i$. Note also that the
forms $\fm$ and $\fm'$ on $\v$ have the same radical $\v^0$.
%So as long as the map $\bar{\;}:\v\rightarrow\vb$ is
%concerned we may work with either of the forms $\fm$ or $\fm'$.

Fix a basis
$\bar{\Pi}_i=\{\bar{\a}_{i1},\ldots,\bar{\a}_{i{\ell_i}}\}$ of
$\bar{R}_i$ and choose a preimage $\a_{ij}\in R_i$ of
$\bar{\a}_{ij}$, $1\leq j\leq {\ell_i}$. Put
$\dot{\Pi}_i=\{\a_{i1},\ldots,\a_{i{\ell_i}}\}$.

Set
$$
\vdi=\span_{\bbbr}\dot{\Pi}_i,\quad
\rdi=\{\dot{\a}\in\vdi\mid\bar{\dot{\a}}\in\bar{R}_i\}\andd
\v^{0}_{i}=\v^0\cap\v_i.
$$
Then
\begin{equation}\label{well2}
\v_i=\vdi\oplus\v^{0}_{i}\andd
\v=(\oplus_{i=1}^{k}\vd_i)\oplus\v^0. \end{equation} (See
(\ref{well}) for this last equality.) Now the map $\bar{\;}$
restricted to $\vd_i$ is an isometry from $\vd_i$ onto $\vb_i$
with respect to the form $\fm'=c_i\fm$ which maps $\rd_i$
bijectively onto $\rb_i$, and so $\rd_i$ is a finite root system
in $\vd_i$ isomorphic to $\rb_i$.

Moreover,
$$
\rtimes_i=\{\da+\d\in R\mid\da\in\rdi^\times,\;\d\in\v^0\}.
$$
Let $\w$ be the Weyl group of $R$ and $\wdot_i$ be the  group
generated by $\{w_{\dot{\a}}\mid\dot{\a}\in\rdi^\times\}$. Since
$(\rdi,\v_j)=\{0\}$ for $i\not=j$ and $(R_i,\v^0)=\{0\}$, by
restriction $\dot{\w}_i$ is isomorphic to the Weyl group of the
finite root system $\rdi$. Moreover, since $R$ contains all
reduced roots of $\rdi$ (see [AABGP, Proposition II.2.11]) and
$w_{r\dot{\a}}=w_{\dot{\a}}$, $\dot{\a}\in\rdi^\times$,
$r\in\bbbr\setminus\{0\}$, $\wdot_i$ is a subgroup of $\w$. It
follows from this that
\begin{equation}\label{z8}
\rtimes_i=\{\da+\d\in R\mid\da\in\rdi^\times,\;\d\in R^0\}.
\end{equation}
In fact if $\da+\d\in R$, $\da\in\rdi^\times$, $\d\in\v^0$, then
$\da=rw_{\dot{\b}_1}\cdots w_{\dot{\b}_m}(\dot{\gamma})$ for some
$\dot{\b}_1,\ldots,\dot{\b}_m,\dot{\gamma}\in\dot{\Pi}_i$ where
$r=2$ or $1$, depending on either $\da/2$ is a root in $\rd_i$ or
not. Since $R$ is $\w$ invariant and isotropic elements are fixed
under the action of the Weyl group we get $\pm r\dot{\gamma}+\d\in
R$. Now consider the $\dot{\gamma}$-string through
$r\dot{\gamma}+\d$ to conclude that $\d\in R$ (see [AABGP,
Proposition II.2.11(b)]). A similar argument shows that if we put
$$R_i^0:=\{\d\in R^0\mid\a+\d\in R\hbox{ for some
}\a\in\rtimes_i\},
$$
then
\begin{equation}\label{2.11b}
\{\d\in R^0\mid\dot{\a}+\d\in R\hbox{ for some
}\dot{\a}\in\rd^{\times}_i\}\sub R^0_i.
\end{equation}
Now let$$R_i=\rtimes_i\cup R_i^0.$$ Then
\begin{equation}\label{w8a}
R=\big(\cup_{i=1}^{k}R_i\big)\cup\riso.
%\andd
%\riso=\emptyset\Longleftrightarrow R=\cup_{i=1}^{k}R_i,
\end{equation}
Since $\v=\span_\bbbr R$, it follows from (\ref{w8a}),
(\ref{well2}) and (\ref{z8}) that
$$\v^0=\span_{\bbbr}R^0.
$$
Set
$$R'_i=R_i\cup(\la R_i\ra\cap R^0).
$$
By (\ref{2.11b}),
\begin{equation}\label{2.11}
(R'_i)^0=\la R_i\ra\cap R^0=\la R^0_i\ra\cap R^0.
\end{equation}

The proof of the following lemma is essentially the same as [ABY,
Lemma 1.2], however for the reader's convenience we give the
details. In what follows we denote by $\la S\ra$, the $\bbbz$-span
of a subset $S$ of a vector space.

\begin{lem}\label{yos}
Let $R$ be a GRRS and $R_1$ be a subset of $R$ with
$\rtimes_1:=R_1\cap\rtimes\not=\emptyset$. Suppose that

(a) $R_1=-R_1$,

(b) $\{\d\in R^0\mid\a'+\d\in R_1\hbox{ for some }\a'\in
R_{1}^{\times}\}\sub R_1$,

(c) $\a'\in R_1$, $\b\in R$, $(\a',\b)\not=0\Longrightarrow \b\in
R_1$.

\noindent Then $R_1$ is a GRRS in its real span. Moreover, if we
set
$$
R_{1}'=R_{1}^{\times}\cup(\la R_1\ra\cap R^0),
$$
then $R_{1}'$ is also a GRRS in the real span of $R_1$.
\end{lem}

\noindent Proof. Since $\rtimes_1\not=\emptyset$, it is enough to
show that axioms (R1)-(R5) hold for $R_1$. Clearly (R1)-(R3) and
(R5) hold for $R_1$. We now check (R4). Let $\a'\in
R_{1}^{\times}$ and $\b'\in R_1$. Since (R4) holds for $R$, it is
enough to show that for $n\in{\bbbz}$,
$$
\b'+n\a'\in R\Longrightarrow\b'+n\a'\in R_1.
$$
Since $\b'\in R_1$, we may assume that $n\not=0$. Assume first
that $n>0$. So let $\b'+n\a'\in R$, $n>0$. If
$\b'+n\a'\in\rtimes$, then $(\b'+n\a',\b')\not=0$ or
$(\b'+n\a',\a')\not=0$. In either cases, we get from (c) that
$\b'+n\a'\in R_1$. Next, let $\b'+n\a'\in R^0$. Since (R4) holds
for $R$ and $n>0$, we have $\b'+(n-1)\a'\in\rtimes$. So repeating
our previous argument we get $\b'+(n-1)\a'\in R_1$. Since
$$
\b'+n\a'+(-\a')=\b'+(n-1)\a'\in R_{1}^{\times}
$$
it follows from (a) and (b) that $\b'+n\a'\in R_1$. If $n<0$ and
$\b'+n\a'\in R$, then $-\b'-n\a'\in R$. Now by the previous step
$-\b'-n\a'\in R_1$ and so by (a), $\b'+n\a'\in R_1$. This
completes the proof of the first assertion.

Next let $R_{1}'$ be as in the statement. Clearly $R_1$ and
$R_{1}'$ have the same real span. Since
$R_{1}^{\times}=(R_{1}')^{\times}$, it is easy to check that
$R_{1}'$ satisfies conditions (a)-(c), and so is a GRRS.
\hfill\qed\vspace{3mm}

\begin{cor}\label{wd}
$R_i$ is an EARS (an irreducible GRRS) and $R'_i$ is an
indecomposable GRRS in $\v_i$.
\end{cor}

\proof It is clear from the way $R_i$ is defined that conditions
(a) and (b) of Lemma \ref{yos} hold for $R_i$. Condition (c) also
holds as $R_i$ is indecomposable.\hfill\qed \vspace{3mm}

We define the {\it rank} of $R_i$ to be the {\it rank} of finite
root system $\rb_i$.
\begin{rem}\label{w8c}
{\rm (i) From [AABGP, Theorem II.2.37] we know that the set of
isotropic roots of an EARS is of the form $S+S$ where $S$ is a
semilattice in the radical of the form. If the nullity $\nu$ is
$1$, then S is a lattice ([AABGP, Corollary II.1.7]). If $\nu=2$
then from [AABGP, II.\S 1] we know that the $\bbbz$-span $\Lam$ of
$S$ is of the form $\Lam=\bbbz\sg_1\oplus\bbbz\sg_2$ where
$\sg_1,\sg_2\in S$ and $\{\sg_1,\sg_2\}$ is a basis of $\v^0$.
Then $S$ is one of the semilattices
$$\Lam,\quad S':=2\Lam\cup(\sg_1+2\Lam)\cup(\sg_2+2\Lam),\quad
S'+\sg_1,\quad S'+\sg_2.$$ Thus $S+S=\Lam$. Therefore if $\nu\leq
2$, the set of isotropic roots is a lattice. Now according to
Corollary \ref{wd} each $R_i$ is an EARS and so $R_i^0$ is a
lattice if $\nu\leq 2$. Thus by (\ref{2.11}), $(R'_i)^0=R_i^0$ and
so $R'_i=R_i$. In particular if $R$ is non-singular, all $R'_i$'s
are also non-singular.

(ii) Even when $R$ is non-singular, the root systems $R'_i$ might
be singular. According to part (i), this only can happen if
$\nu\geq 3$. To see an example let $\nu=3$,
$\v={\bbbr}\a_1\oplus{\bbbr}\a_2\oplus{\bbbr}
\sg_1\oplus{\bbbr}\sg_2\oplus{\bbbr}\sg_3$. Define a positive
semi-definite bilinear form on $\v$ by letting $\sg_i$'s to be
isotropic and $(\a_i,\a_j)=2\d_{ij}$. Set
$$
R_1=(S+S)\cup (\pm\a_1+S)\andd R_2=\Lam\cup (\pm\a_2+\Lam),
$$
where
$$\Lam={\bbbz}\sg_1\oplus{\bbbz}\sg_2\oplus{\bbbz}\sg_3\andd
S=2\Lam\cup(\sg_1+2\Lam)\cup(\sg_2+2\Lam)\cup(\sg_3+2\Lam).
$$
By [AABGP, Theorem II.2.37], $R_1$ and $R_2$ are  GRRS of type
$A_1$ in the real span of $R_1$ and $R_2$, respectively. Set
$R=R_1\cup R_2$. Then $R$ is a non-singular GRRS in $\v$. However
$R'_1=\Lam\cup(\pm\a_1+S)$ is a singular GRRS.}
\end{rem}

We now return to the GRLA  $\gg$ and the corresponding GRRS $R$.
Set
$$\begin{array}{l}
\hh(R)=\span_{\bbbc}\{t_\a\mid\a\in R\},\vspace{3mm}\\
\hh^0(R)=\{h\in\hh(R)\mid \big(h,\hh(R)\big)=\{0\}\}\hbox{ and}\vspace{3mm}\\
\hdi=\span_{\bbbc}\{t_\a\mid\a\in\dot{R}_i\}=
\span_{\bbbc}\{t_\a\mid\a\in\rdi^\times\}.
\end{array}
$$

\begin{lem}\label{w7}
(i) The form $\fm$ on $\hh$ restricted to $\hdi$ is nondegenerate.

(ii) $\hh(R)=\hh^0(R)\oplus\big(\bigoplus_{i=1}^{k}{}\hdi\big)$.

(iii) $\hh^0(R)=\span_{\bbbc}\{t_\d\mid\d\in R^0\}.$
\end{lem}

\proof (i) Let $h=\sum (k_j+ik'_j)t_{\a_j}\in\hd_t$ be in the
radical of the form $\fm$ restricted to $\hd_t$, where
$\a_j\in\rd_t$, $k_j,k'_j\in\bbbr$ and $i$ is the complex number
with $i^2=-1$. In particular
$$(\sum k_j\a_j,\vd_t)+i(\sum k'_j\a_j,\vd_t)=\{0\}.$$
Multiplying both sides by the scalar $c_t$ (see (\ref{scaled})) we
obtain $$(\sum k_j\a_j,\vd_t)_t+i(\sum k'_j\a_j,\vd_t)_t=\{0\}.$$
Since $\fm_t$ is real valued and positive definite (in particular
nondegenerate) on $\vd_t$ we get $\sum k_j\a_j=0$, $\sum
k'_j\a_j=0$ and so $h=0$.

(ii) Since $\rdi$ and $\dot{R}_j$ are orthogonal if $i\not=j$, we
have $(\hdi,\dot{\hh}_j)=\{0\}$. So by part (i),
$\hdi\cap(\sum_{j\not=i}\dot{\hh}_j)=\{0\}$. In particular, the
sum $\sum_{i=1}^{k}\dot{\hh}_i$ is direct and the form restricted
to $\oplus_{i=1}^{k}\dot{\hh}_i$ is nondegenerate. Thus
$(\bigoplus_{i=1}^{k}\hdi)\cap\hh^0(R)=\{0\}$. Note that if $\a\in
R$, then $\a=\dot{\a}+\d$, where $\dot{\a}\in\rdi$ and
$\d\in\v^0$, for some $1\leq i\leq k$. Now $t_{\dot{\a}}\in\hdi$
and $t_\d\in\hh^0(R)$. So $t_\a\in\hdi\bigoplus\hh^0(R)$.

(iii) Let $\a\in R$. By (\ref{z8}), $\a=\dot{\a}+\d$ where
$\dot{\a}\in\rdi$ for some $1\leq i\leq k$ and $\d\in R^0$.
Therefore $t_\a$ is in the span of
$\{t_{\dot{\a}},t_\d\mid\dot{\a}\in\dot{{\Pi}}_i,\d\in R^0\}$. It
follows that if $h\in\hh^0(R)\sub\hh(R)$, then $h=\dot{h}+h^0$
where $\dot{h}\in\bigoplus\hdi$ and
$h^0\in\span_{\bbbc}\{t_\d\mid\d\in R^0\}$. But by part (ii),
$\dot{h}=0$ and $h=h^0$.\hfill\qed\vspace{3mm}

Set
$$\hi^0:=\span_{\bbbc}\{t_\d\mid\d\in R^0,\;\da+\d\in
R\hbox{ for some }\da\in\rdi^\times\}\subset\hh(R)^0.
$$
From (\ref{2.11b}) it follows that
$\hh^0_i\sub\span_{\bbbc}\{t_\d\mid\d\in R^0_i\}.$ By Corollary
\ref{wd}, $R_i$ is an EARS. So by [AABGP, Corollary II.2.31] if
$\d\in R^0_i$, then $\d=\d_1+\d_2$, where
$t_{\d_1},t_{\d_2}\in\hh^0_i$. Thus $t_\d\in\hh^0_i$ and
\begin{equation}\label{2.11c}
\hh^0_i=\span_{\bbbc}\{t_\d\mid\d\in
R^0_i\}=\span_{\bbbc}\{t_\d\mid\d\in (R'_i)^0\}.
\end{equation}

Define
$$
\hd=\bigoplus_{i=1}^{k}\hdi\andd\hh^0=\sum_{i=1}^{k}\hi^0\subseteq\hh(R)^0.
$$

Recall from Definition \ref{qssla} that the subalgebra $\gc$ of
$\gg$ is the subalgebra of $\gg$ generated by root spaces $\ga$,
$\a\in\rtimes$.
\begin{lem}\label{w8} $\hd\oplus\hh^0=
\sum_{\a\in\rtimes}[\ga,\gg_{-\a}]=\gc\cap\hh$.
\end{lem}

\proof It follows immediately from the definition of $\gc$ and
(\ref{w1}) that the second equality holds. Since $\gg$ satisfies
axioms GR1-GR4 of an EALA, we have from (\ref{w4}) that if
$\a\in\rtimes$ then $[\ga,\gg_{-\a}]={\bbbc}t_{\a}$. So if
$\da\in\dot{\Pi}_i\subset\rtimes_i\subset\rtimes$, then
$t_{\dot{\a}}\in[\gg_{\dot{\a}},\gg_{-\dot{\a}}]$. Thus
$t_{\da}\in \sum_{\a\in\rtimes}[\ga,\gg_{-\a}]$ for all
$\da\in\rd_i$ and so $\hdi\sub\sum_{\a\in\rtimes}[\ga,\gg_{-\a}]$.
Next let $\d\in R^0_i$. Then $\a+\d\in R$ for some
$\a\in\rtimes_i$. Also by (\ref{z8}) and (\ref{2.11b})
$\a=\dot{\a}+\eta$ for some $\dot{\a}\in\rdi\setminus\{0\}$ and
$\eta\in R^0_i$. Now $t_{\dot{\a}},t_\a,
t_{\a+\d}\in\sum_{\a\in\rtimes}[\ga,\gg_{-\a}]$. Thus the sum
contains $t_\d$ too, and so contains $\hd\oplus\hh^0$.

Conversely, let $\a\in\rtimes$, then $\a\in R_i$ for some $i$. By
(\ref{z8}) and (\ref{2.11b}) $\a=\dot{\a}+\d$ where
$\dot{\a}\in\rdi\setminus\{0\}$ and $\d\in R^0_i$. Then
$t_\a=t_{\dot{\a}}+t_\d\in\hd\oplus\hh^0_i$.
\hfill\qed\vspace{3mm}

It can be read from proof of Lemma \ref{w8} that
\begin{cor}\label{w8b} $\hdi\oplus\hi^0=\sum_{\a\in \rtimes_i}[\ga,\gg_{-\a}]$.
\end{cor}

Set
\begin{equation}\label{w15}
R^0_0=\riso\setminus
(\cup_{i=1}^{k}R'_i)=\riso\setminus\big(\cup_{i=1}^{k}(R'_i)_{\hbox{iso}}\big).
\end{equation}
%and
%$$
%\hh^0_0=\span_{\bbbc}\{t_\d\mid\d\in R^0_0\}.
%$$

%\begin{lem}\label{w10a}
%$\hh^0(R)=\hh^0\oplus\hh^0_0$.
%\end{lem}

%\proof It is enough to show the inclusion $\sub$. By Lemma
%\ref{w7}(iii), we must show that $t_\d\in\hh^0\oplus\hh^0_0$ for
%any $\d\in R^0$. If $\d\in\rniso$, then $\d+\a\in\rtimes_i$ for
%some $i$ and some $\a\in\rtimes_i$, so $t_\d\in\hh^0_i\sub\hh^0$.
%If $\d\in\riso\setminus (\cup_{i=1}^{k}R'_i)$, then
%$t_{\d}\in\hh^0_0$. If $\d\in\riso\cap R'_i$ for some $i$, then
%$\d\in\la R^0_i\ra$, but $t_\gamma\in\hh^0$ for any $\gamma\in
%R^0_i$ and so $t_\d\in\hh^0$.\qed

Since $\fm$ is nondegenerate on $\hh$ and $\hd$, we have
$\hh=\hd\oplus\hd^\perp$ where $\hd^\perp$ is the orthogonal
complement of $\hd$ in $\hh$. Since
$(\hd\oplus\hh^0,\hh^0)=\{0\}$, there exists a subspace $\dd$ of
$\hd^\perp$ such that
\begin{equation}\label{www}
\begin{array}{l}
\dim\dd=\dim\hh^0,\vspace{3mm}\\
(\dd,\dd\oplus\hd)=\{0\},\vspace{3mm}\\
\fm\hbox{ is nondegenerate on
}\hd\oplus\hh^0\oplus\dd\vspace{3mm}.\\
%(\dd,\hh^0_0)=\{0\}.
\end{array}
\end{equation}
%To see this write $\hh=\hd\oplus\hh^0\oplus\hh^0_0\oplus M$ for
%some vector space $M$. Consider a basis $d'_1,\ldots,d'_t$ of
%$\hh^0$ and define $x_i\in\hh^\star$ by $x_i(\hd+\hh^0_0+M)=0$ and
%$x_i(d'_j)=\d_{ij}$. Let $d_i$ be the unique element in $\hh$
%representing $x_i$ and let $\dd=\oplus\bbbc d_i$. Then
%$(\dd,\hd+\hh^0_0+M)=0$. Note that since $(\hh^0,\hh^0)=0$, we may
%replace $d_i$ with $d_i-\frac{(d_i,d_i)}{2}d'_i$ to assume that
%$(d_i,d_i)=0$. Next set $h_1=d_1$, $h_2=d_2-(d_1,d_2)d'_2$ and in
%general $h_j=d_j-\sum_{i=1}^{j-1}(d_j,d_i)d'_i$. Now replace
%$\dd$ with the new space $\oplus\bbbc h_i$. Then $\dd$ satisfies
%all the requirements in (\ref{www}).
%In particular, $\hh^0_0\cap \dd=\{0\}$.
Next, let $\w$ be the orthogonal complement of
$\hd\oplus\hh^0\oplus\dd$ in $\hh$, then we have
\begin{equation}\label{w11}
\begin{array}{l}
\hh=\hd\oplus\hh^0\oplus\dd\oplus\w\vspace{3mm}\\
\big(\hd\oplus\hh^0\oplus\dd,\w\big)=\{0\}\vspace{3mm}\\
%\big(\hh^0_0\oplus\w_0^0,\w_0)=\{0\}\vspace{3mm}\\
\fm\hbox{ is nondegenerate
 on } \w.
\end{array}
\end{equation}
Now consider a basis $B=\{h_1,\ldots,h_m\}$ of $\hh^0$ such that
$B$ contains a basis of $\hh^0_i$, for each $i$. Using
(\ref{www}), we may pick a basis $B'=\{d_1,\ldots,d_m\}$ of $\dd$
such that $(h_i,d_j)=\d_{ij}$. Let
$$\dd_i=\span_{\bbbc}\{d_j\in B'\mid h_j\in\hh^0_i\}.$$
Then $\dd=\sum_{i=1}^{k}\dd_i$ and
\begin{equation}\label{w10}
\begin{array}{l}
\dim\dd_i=\dim\hh^0_i,\vspace{3mm}\\
(\dd_i,\dd\oplus\hd)=\{0\},\vspace{3mm}\\
\fm\hbox{ is nondegenerate on }\hd\oplus\hh^0_i\oplus\dd_i.
\end{array}
\end{equation}
%It follows that $\dd=\sum_{i=1}^{k}\dd_i$. It is enough to show
%that $\dim\sum\dd_i=\dim\hh^0$. For this pick a basis
%$h_{11},h_{12},\ldots,h_{1m_{1}},
%\ldots,h_{k1},h_{k2},\ldots,h_{km_{k}}$ of $\hh^0$ where
%$h_{ji}\in\hh^0_j$ and $\{h_{ij}\}_{i,j}$ contains a basis of
%$\hh^0_i$ for each $i$. Define $\a_{ij}\in
%(\hd\oplus\hh^0\oplus\dd)^\star$ by
%$\a_{ij}(h_{kl})=\d_{ik}\d_{jl}$, $\a_{ij}(\dd)=\{0\}$ and
%$\a_{ij}(\hd)=\{0\}$. We may consider $\a_{ij}$ as an element of
%$(\hd\oplus\hh^0_i\oplus\dd_i)^\star$ too. Let $d_{ij}$ be the
%unique element in $\hd+\hh^0+\dd_i$ which represent $\a_{ij}$ on
%$\hd\oplus\hh^0_i\oplus\dd_i$. It follows that $d_{ij}\in\dd_i$.
%Then $d_{ij}$'s are linearly independent. In fact, if $\sum
%m_{ij}d_{ij}=0$, then $\a_{lk}(\sum m_{ij}d_{ij})=0$. This gives
%$m_{kl}=0$.
Note that if $\a\in R'_i$, then $\a=\da+\d$ where
$t_{\da}\in\hd_i$ and $t_\d\in\hh^0_i$. Since
$(\hd_i\oplus\hh^0_i,\hd_j)=\{0\}$ for $i\not=j$, we have
$\a(\hd)=\da(\hd)=\da(\hd_i)$. Now from (\ref{www}), (\ref{w11}),
(\ref{w10}) and the way the spaces $\dd_i$'s are defined, we have
$\d(\hh)=\d(\dd)=\d(\dd_i)$. Thus
\begin{equation}\label{rev1}
\a(\hh)=\a(\hh_i)=\a(\hd_i\oplus\dd_i)\quad\hbox{for}\quad\a\in
R'_i.
\end{equation}
It follows from this that for $\a\in R'_i\setminus\{0\}$,
\begin{equation}\label{rev2}
\gg_\a=\{x\in\sum_{\b\in
R'_i\setminus\{0\}}\gg_\b\mid[h,x]=\a(h)x\hbox{ for all
}h\in\hh^i\}.
\end{equation}
In fact, it is clear that $\gg_\a$ is a subset of the right hand
side. To see the reverse inclusion let
$x=x_{\b_1}+\cdots+x_{\b_t}$ be an element of the right hand side,
where $\b'_j\in R'_i\setminus\{0\}$ for all $j$, $x_\b\in\gg_\b$.
Let $h=\dot{h}+h^0+d+h^0_0+w^0_0+w_0$ be an arbitrary element of
$\hh$ in the form (\ref{w11}), where
$\dot{h}=\sum_{j=1}^{k}\dot{h}_j$, $\dot{h}_j\in\hd_j$. With
respect to the basis $B'$ of $\dd$ we may write $d=d'+d''$ where
$d'\in\dd_i$ and $d''\in\span_{\bbbc}\{d_j\mid d_j\not\in
\dd_i\}$. Now using (\ref{rev1}) we have
$$
[h,x]=\sum_{j=1}^{t}\b_j(h)x_{\b_j}=
\sum_{j=1}^{t}\b_j(\dot{h_i}+d')x_{\b_j}
=[\dot{h_i}+d',x]=\a(\dot{h_i}+d')x=\a(h)x.
$$
So $x\in\gg_\a$.

Starting from each $R_i$, we now would like to construct a
generalized reductive subalgebra $\gg^i$ of $\gg$ which is
indecomposable. For this set
\begin{equation}\label{w13a}
\hi=\hdi\oplus\hi^0\oplus\dd_i\andd \gg^i=\hi\oplus\sum_{\a\in
R'_i\setminus\{0\}}\ga \end{equation}

\begin{pro}\label{w13b}
(i) $(\gg^i,\fm,\hi)$ is an indecomposable generalized reductive
subalgebra of $\gg$.

(ii) $\gg^i_c\cap\hh_i=\hdi\oplus\hh^0_i$.
\end{pro}

\proof First we must show that $\gg^i$ is a subalgebra of $\gg$.
Note that $R'_i$, $R_i$ and $\rtimes_i$ have the same linear span.
It then follows from (\ref{w4}) and Corollary \ref{w8b} that
$$
[\gg^i,\gg^i]\cap\hh=\sum_{\a\in R'_i}[\ga,\gg_{-\a}]=\sum_{\a\in
\rtimes_i}[\ga,\gg_{-\a}]=\hdi\oplus\hi^0\subset\gg^i.
$$
Since $\gg_0=\hh$ acts diagonally on $\gg$, we have
$[\hh,\gg^i]\subset\gg^i$. Thus it only remains to show that if
$\a,\b\in R'_i\setminus\{0\}$ and $\a+\b\in R\setminus\{0\}$, then
$\a+\b\in R'_i$. If $\a+\b$ is isotropic, it is clear from
definition of $R'_i$ that $\a+\b\in R'_i$. If $\a+\b$ is
non-isotropic, then it can not be orthogonal to both $\a$ and $\b$
and so $\a+\b\in R_i\subset R'_i$, since $R_i$ is indecomposable.
Thus $\gg^i$ is a subalgebra of $\gg$. Moreover, from (\ref{rev1})
and (\ref{rev2}) we have $\hh_i=C_{\gg^i}(\hh_i)$ and
$$\gg_\a=\{x\in\gg^i\mid [h,x]=\a(h)x\hbox{ for all }h\in\hh^i\}
\quad(\a\in R'_i\setminus\{0\}).
$$
Thus
\begin{equation}\label{w14}
\gg^i=\sum_{\a\in R'_i}(\gg^i)_\a,\hbox{ where
}(\gg^i)_0=\hi\hbox{ and }(\gg^i)_{\a}=\ga\hbox{ for }\a\not=0.
\end{equation}
Next we must show that GR1-GR5 and GR6a hold for $\gg^i$. By
(\ref{w10}), (\ref{w2}) and Corollary \ref{wd} GR1 holds for
$\gg^i$. Considering $R'_i$ as a subset of $\hi^\star$, we see
from (\ref{w14}) that elements from $\hi$ act diagonally on
$\gg^i$ via the adjoint representation. So GR2 holds as
$C_{\gg^i}(\hi)=\hi$. Validity of GR3 and GR4 for $\gg^i$ follows
from the fact that these axioms hold for $\gg$ and that
$R'_i\subset R$, and it is clear that GR5 holds. Finally, GR6a
holds, since $R'_i$ is indecomposable. Part (ii) follows from
Corollary \ref{w8b}. \hfill\qed\vspace{3mm}

Put
\begin{equation}\label{w16}
\i=\sum_{\d\in R^0_0}\gg_{\d}.
\end{equation}
From (\ref{w1}), (\ref{w8a}), (\ref{w11}) and (\ref{w16}) we have
\begin{equation}\label{w17}
\gg=\sum_{i=1}^{k}\gg^i\oplus\w\oplus\i.
\end{equation}
Here $\w$ and each $\gg^i$'s are Lie subalgebras of $\gg$ and $\i$
is a subspace of $\gg$. The direct sums appearing in (\ref{w17})
are just sums of vector spaces. We now would like to investigate
the Lie bracket between these spaces, at the level of cores.

\begin{lem}\label{w19}
(i) $\gc=\sum_{i=1}^{k}\gg^i_c$ where $[\gg^i_c,\gg^j_c]=\{0\}$
for $i\not=j$. In particular, $\gg^i_c$ is an ideal of $\gg$, for
each $i$.

(ii) If $x=\sum_{i=1}^{k}x_i\in{\mathcal Z}(\gc)$ where
$x_i\in\gc^i$, then $x_i\in{\mathcal Z}(\gc^i)$ for each $i$. In
particular, ${\mathcal Z}(\gc)=\sum_{i=1}^{k}{\mathcal
Z}(\gg^i_c)$ and as Lie algebras
$$\frac{\gc}{{\mathcal
Z}(\gc)}\cong\bigoplus_{i=1}^{k}\frac{\gg^i_c}{{\mathcal
Z}(\gg^i_c)}.$$
\end{lem}

\proof (i) Let $i\not=j$, $\a\in\rtimes_i$ and $\b\in\rtimes_j$.
Then $(\a,\b)=0$ and so $\a+\b$ is orthogonal to neither of $\a$
and $\b$. Thus $\a+\b$ is not a root of $\gg$. This shows that
$[\gc^i,\gc^j]=\{0\}$ for $i\not=j$, and that $\gc^i$ is an ideal
of $\gc$. Since $\gc^i$ is perfect, it follows from the Jacobi
identity that $\gc^i$ is an ideal of $\gg$. In particular,
$\sum_{i=1}^{k}\gg^i_c$ is a subalgebra of $\gg$ containing all
non-isotropic root spaces. Clearly any subalgebra of $\gg$
containing all non-isotropic root spaces must contain this sum.
Thus $\gc=\sum\gg^i_c$.

(ii) Let $x$ be as in the statement. Then by part (i) for each
$i$, we have
$$ \{0\}=[x,\gg_c^i]=\sum_{j=1}^{k}[x_j,\gg_c^i]=[x_i,\gg_c^i].$$
So $x_i\in{\mathcal Z}(\gg_c^i)$. It now follows that
$\bigoplus_{i=1}^{k}{\mathcal Z}(\gg_c^i)$ is the kernel of the
epimorphism
$$\begin{array}{c}
\displaystyle{\bigoplus_{i=1}^{k}\gg_c^i\longrightarrow\frac{\gg_c}{{\mathcal
Z}(\gg_c)}}\vspace{2mm}\\
\displaystyle{(x_i)\longmapsto\sum x_i+{\mathcal
Z}(\gg_c)}.\end{array}
$$
\hfill\qed\vspace{3mm}

\begin{lem}\label{w20} (i) $\w\subseteq
C_{\gg}(\sum_{i=1}^{k}\gg^i)$.

(ii) $\sum_{\d\in\riso}\gg_\d\subseteq C_{\gg}(\gc)$. In
particular, $\w\oplus\i\subseteq C_{\gg}(\gc)$.
\end{lem}

\proof (i) Let $\a\in R'_i$. We have $\la
R'_i\ra=\la\rtimes_i\ra$. It follows from this that $t_\a\in\hi$.
Then from (\ref{w11}) we have
$$[\w,\ga]=\a(\w)\ga=(t_\a,\w)\ga\subseteq
(\hi,\w)\ga=\{0\}.
$$
%A similar argument shows $[\w^0_0,\ga]=\{0\}$.

(ii) The first part of the statement follows from  Proposition
\ref{w6}(ii). The second part of the statement holds now by part
(i).
%(iii) We have $[\w,\w^0_0]\subseteq [\hh,\hh]=\{0\}$. If
%$\d\in\riso\setminus (\cup_{i=1}^{k}R'_i)$, then $t_\d\in\hh^0_0$.
%So from (\ref{w11}) we have
%$$[\w,\gg_\d]=\d(\w)\gg_\d=(t_\d,\w)\gg_\d\subseteq
%(\hh^0_0,\w)\gg_\d=\{0\}.
%$$
\hfill\qed\vspace{3mm}

Let us summarize the results obtained in the following theorem.

\begin{thm}\label{main} Let $(\gg,\fm,\hh)$ be a GRLA with
corresponding root system $R$. Then

(i) $R=(\cup_{i=1}^{k}R_i)\cup\riso$ where for each $i$, $R_i$ is
an EARS. Moreover $R'_i=R_i\cup(\la R_i\ra\cap R^0)$ is an
indecomposable GRRS.

(ii) For $1\leq i\leq k$, there exists a subspace $\dd_i$ of $\hh$
such that if
$\hh_i:=\dd_i\oplus\sum_{\a\in\rtimes_i}[\ga,\gg_{-\a}]$, and
$\gg^i:=\hh_i\oplus\sum_{\a\in R'_i}\gg_{\a}$, then
$(\gg^i,\fm,\hh_i)$ is an indecomposable generalized reductive
subalgebra of $\gg$. Moreover,
$$\dim\dd_i=\dim(\sum_{\a\in\rtimes_i}
[\ga,\gg_{-\a}])-\hbox{rank}({R}_i).$$

(iii) $\hh$ has a decomposition as in (\ref{w11}). Moreover
$\gg=\sum_{i=1}^{k}\gg^i\oplus\w\oplus\i$, where $\w$ and $\i$ are
subspaces of $\gg$ defined by (\ref{w16}). Moreover,
$[(\sum_{i=1}^{k}\gg^i)\oplus\w,\w]=\{0\}$ and $[\i,\gc]=\{0\}$.

(iv) If $i\not=j$, then $[\gg^i_c,\gg^j_c]=\{0\}$ and
$\gc=\sum_{i=1}^{k}\gg^i_c$. In particular, $\gg^i_c$ is an ideal
of $\gg$. Moreover, ${\mathcal Z}(\gc)=\sum_{i=1}^{k}{\mathcal
Z}(\gg^i_c)$ and
$$\frac{\gc}{{\mathcal
Z}(\gc)}\cong\bigoplus_{i=1}^{k}\frac{\gg^i_c}{{\mathcal
Z}(\gg^i_c)}.$$

(v) If $\pi$ is the projection map $\v\rightarrow\v^0$ (with
respect to the decomposition (\ref{well2})), then for $i\not= j$,
$$
[\gg^i,\gg^j]\sub(\hh^0_i\cap\hh^0_j)\oplus\sum_{\a\in R'_i\cup
R'_j,\;\pi(\a)\not=0}\ga+\sum_{\a',\b'}[\gg_{\a'},\gg_{\b'}],
$$
where $(\a',\b')\in (\rtimes_i\times(R'_j)^0\setminus\{0\})\cup
(\rtimes_j\times(R'_i)^0\setminus\{0\})$. In particular,
$[\gg^i,\gg^j]=\{0\}$ if $\v^0=\{0\}$.

(vi) If $\gg$ is tame, then $\w=\{0\}$ and $\i=\{0\}$.
\end{thm}

\proof (i) See (\ref{w8a}) and Corollary \ref{wd}. (ii) See
Corollary \ref{w8b}, (\ref{w10}), (\ref{w13a}) and Proposition
\ref{w13b}. (iii) See (\ref{w17}) and Lemma \ref{w20}. (iv) See
Lemma \ref{w19}.

(v) We must check $[(\gg^i)_\a,(\gg^j)_\b]$ for $\a\in R'_i,\b\in
R'_j$. First let $\a=0$ and $\b\in R'_j\setminus\{0\}$. Then
$(\gg^i)_\a=\hh_i$ and $(\gg^j)_\b=\gg_\b$. Since
$t_{\b-\pi(\b)}\in\hd_j$, we have from (\ref{www}) that
$$[\hh_i,\gg_{\b-\pi(\b)}]=(t_{\b-\pi(\b)},\hh_i)\subset
(\hd_j,\hh_i)=\{0\}.
$$
Thus
$$
[\hh_i,\sum_{\b\in R'_j\setminus\{0\}}\gg_\b]\sub\sum_{\b\in
R'_j,\;\pi(\b)\not=0}\gg_\b.
$$
 Next let $\a\in R'_i\setminus\{0\}$, $\b\in R'_j\setminus\{0\}$.
 If $-\b=\a=\d\in (R'_i)^0\cap (R'_j)^0$, then by
(\ref{w4}) and (\ref{2.11c}), we have
$[(\gg^i)_\a,(\gg^j)_\b]=[\gg_\d,\gg_{-\d}]={\bbbc}t_\d\sub\hh^0_i\cap\hh^0_j$.
It then follows from (\ref{w19})(i) that for $0\not=\a\in R'_i$
and $0\not=\b\in R'_j$,
$$
[\ga,\gg_\b]\sub\hh^0_i\cap\hh^0_j\oplus\sum_{\a',\b'}[\gg_{\a'},\gg_{\b'}],
$$ where $\a',\b'$ are as in the statement.

(vi) If $\gg$ is tame then $C_{\gg}(\gc)\sub\gc$ and by [ABP,
Lemma 3.62], $\riso=\emptyset$ (and so $\i=\{0\}$). By part (iii),
$\w\sub C_{\gg}(\gc)\cap\hh\sub\gc\cap\hh$. But by Lemma \ref{w8},
$\gc\cap\hh=\hd\oplus\hh^0$. So from (\ref{w11}) we get
$\w=\{0\}$. \hfill\qed\vspace{3mm}

\begin{cor}\label{main2}
Let $\gfh$ be a non-singular GRLA of nullity $\nu$ with root
system $R$. Then $\gg=\sum_{i=1}^{k}\gg^i\oplus\w$, where
$\gg^i$'s are indecomposable GRLAs with $[\gg^i_c,\gg^j_c]=\{0\}$,
for $i\not=j$, and $\w$ is contained in the centralizer of $\gg$.
Moreover,

(i) if $\nu\leq 2$, then each $\gg^i$ is an EALA,

(ii) if $\nu=0$, then $\gg$ is a finite dimensional reductive Lie
algebra.
\end{cor}

\proof  By assumption $\riso=\emptyset$ and so $\i=\{0\}$. If
$\nu\leq 2$, we have from Remark \ref{w8c}(i) that $R'_i=R_i$.
That is $(R'_i)_{\iso}=\emptyset$. Thus $\gg^i$ satisfies GR6b.

If $\nu=0$, then from Theorem \ref{main}(v), we have that the sum
in the statement is direct, that is
$\gg=(\bigoplus_{i=1}^{k}\gg^i)\oplus\w$. Moreover, $R_i=\rd_i$ is
an irreducible finite root system and
$$
\gg^i=\bigoplus_{\a\in\rd_i}(\gg^i)_\a\hbox{ where
}(\gg^i)_0=\hh^i\hbox{ and }(\gg^i)_\a=\gg_\a\hbox{ for }\a\not=0.
$$
Since $\hh^0_i=\{0\}=\dd_i$, we have
$\dim\hh_i=\dim\hd_i=\hbox{rank}\rd_i$. Now it follows from
Theorem \ref{aabgp} and Serre's Theorem that each $\gg^i$ is a
finite dimensional simple Lie algebra over the field of complex
numbers. That is $\gg$ is a reductive Lie algebra.
\hfill\qed\vspace{3mm}

It is worth mentioning that the basic structural properties of an
EALA essentially come from its core modulo its center (see [AG,
Proposition 1.28], [A2, Proposition 1.28] and [N2, Theorem 6]).
Therefore Theorem \ref{main}(iv) together with Corollary
\ref{main2}(i) suggest that the structural properties of a
generalized reductive Lie algebra $\gg$ can be obtained from the
indecomposable subalgebras $\gg^i$.

%\begin{exa}\label{exa2}
%For $1\leq i\leq k$ let $(\gg^i,\fm,\hh_i)$ be an EALA
%(irreducible GRLA) and let $\w$ be an abelian Lie algebra. Set
%$\gg=\oplus_{i=1}^{k}\gg^i\oplus\w$. Consider a nondegenerate
%symmetric bilinear form on $\w$. Putting together the forms on
%each piece we get a nondegenerate symmetric invariant form on
%$\gg$. Set $\hh=\oplus_{i=1}^{k}\hh_i\oplus\w$. Then $\gfh$ is a
%nonsingular GRLA.
%\end{exa}

\section{\bf On the classification of GRLAs}
\setcounter{equation}{0} In this section, we show that the core
modulo center of an indecomposable GRLA is a centerless Lie torus.
Therefore by Theorem \ref{main}(ii),(iv) the core modulo center of
a GRLA is a direct sum of centerless Lie tori. For the
classification of centerless Lie tori of types $A_1$ and $A_2$ see
[Y] and [BGKN] respectively. For simply laced types of rank $\geq
3$ see [BGK]. For types $B_\ell$, $C_\ell$, $F_4$ and $G_2$ see
[AG]. For type $BC_{\ell}$ ($\ell\geq 3$) see [ABG]. Finally, for
type $BC_1$ see [AFY]. The classification for type $BC_2$ is open.

Let us recall the definition of a Lie torus over $\bbbc$,
introduced in [N1]. Let $\ll$ be a complex Lie algebra, $\rd$ be
an irreducible finite root system and $\Lam$ be a free abelian
group of finite rank. Denote the set of indivisible roots of $\rd$
by $\rd_{\ind}$, that is
$\rd_{\ind}=\{\da\in\rd\mid\frac{1}{2}\da\not\in\rd\}.$ Then, the
Lie algebra $\ll$ is called a {\it Lie torus of type} $(\rd,\Lam)$
if it satisfies the following axioms:

\noindent (LT1) $\ll$ has a $(\la\rd\ra\oplus\Lam)$-grading of the
form
$$
\ll=\bigoplus_{\da\in\lrd,\lam\in\Lam}\ll^{\lam}_{\da},\quad[\ll_{\da}^{\lam},
\ll^{\mu}_{\db}]\sub\ll^{\lam+\mu}_{\da+\db},\quad
\hbox{satisfying }\ll^{\lam}_{\da}=\{0\}\hbox{ if }\da\not\in\rd.
$$

\noindent (LT2) For $\da\in\rd^\times:=\rd\setminus\{0\}$ and
$\lam\in\Lam$, we have

(i) $\dim\ll^{\lam}_{\da}\leq 1$, with $\dim \ll^{0}_{\da}=1$ if
$\da\in\rd_{\ind}$,

(ii) if $\dim L_{\da}^{\lam}=1$, then there exists
$(e^{\lam}_{\da},f^{\lam}_{\da})\in\ll^{\lam}_{\da}\times\ll^{-\lam}_{-\da}$
such that
$h^{\lam}_{\da}:=[e^{\lam}_{\da},f^{\lam}_{\da}]\in\ll^{0}_{0}$
acts on $x\in\ll^{\mu}_{\db}$ $(\db\in\rd,\mu\in\Lam)$ by
$$
[h^{\lam}_{\da},x]=(\db,{\da}^\vee)x,
$$
where $(\db,{\da}^\vee)$ is the Cartan integer of $\db,\da$.

\noindent (LT3) For $\lam\in\Lam$ we have
$\ll^{\lam}_{0}=\sum_{\da\in\rd^\times,\mu\in\Lam}[\ll^{\mu}_{\da},\ll^{\lam-\mu}_{-\da}].$

\noindent (LT4)
$\Lam=\la\{\lam\in\Lam\mid\ll^{\lam}_{\da}\not=\{0\}\hbox{ for
some }\da\in\rd\}\ra.$

 We start with an indecomposable GRLA
$(\gg,\fm,\hh)$, that is $\gg$ satisfies axioms GR1-GR5 and GR6a.
Let $R$ be the root system of $\gg$. According to Section 1, we
have
$$
R=R_t\cup\riso,\quad\hbox{where}\quad R_t=\rtimes\cup R_t^0,
$$
with $$R_t^0=\{\d\in R^0\mid\d+\a\in R\hbox{ for some
}\a\in\rtimes\}.
$$
Moreover, $R_t$ is an EARS (an irreducible GRRS) and $\Lam_t=\la
R^0_t\ra$ is a lattice. So there exists an irreducible finite root
system $\rd$ with
$$
R_t\sub\rd+\Lam_t\andd \rd_{\ind}\sub\rtimes.
$$
(See [AABGP, Proposition II.2.11].) Thus
$$
R\sub\rd+\Lam_t+\Lam_0,\quad\hbox{where}\quad\Lam_0=\la\riso\ra.
$$
Set
$$
\Lam=\la R^0\ra=\Lam_t+\Lam_0=\la R_t^0+\riso\ra.
$$
So $\la R\ra=\la\rd\ra\oplus\Lam$.

\begin{thm}\label{classification}
Let $(\gg,\fm,\hh)$ be an indecomposable GRLA. Then the Lie
algebra $\gc/{\mathcal Z}(\gc)$ is a Lie torus of type
$(\rd,\Lam_t)$.
\end{thm}

\proof We start by proving that axioms LT1 and LT2 of a Lie torus
hold for $\gg$ with respect to a grading based on the abelian
group $\lrd\oplus\Lam$. From Section 1, we know that
$$
\gg=\bigoplus_{\a\in R}\gg_\a=\bigoplus_{\da\in\lrd,\lam\in\Lam}
\gg_{\da+\lam}
$$
is a $(\lrd\oplus\Lam)$-grading for $\gg$ with
$\gg_{\da+\lam}=\{0\}$ if $\da\not\in\rd.$ Thus by considering
$\ll_{\da}^{\lam}:=\gg_{\da+\lam}$, $(\da\in\lrd,\lam\in\Lam)$, we
see that axiom LT1 of a Lie torus holds for $\gg$ with respect to
this grading. Next let $\da\in\rd^\times$ and $\lam\in\Lam$. If
$\da+\lam\in R$, then $\da+\lam\in\rtimes$ and so by Theorem
\ref{aabgp}(d), $\dim\gg_{\da+\lam}=1$. If $\da+\lam\not\in R$,
then $\dim\gg_{\da+\lam}=0$. Moreover, $\rd_{\ind}\sub R$ and so
$\dim\gg_{\da}=1$ for $\da\in\rd_{\ind}\setminus\{0\}$. Thus part
(i) of the axiom LT2 holds for $\gg$.

Next note that if $\da\in\rd^\times$, $\lam\in\Lam$ and
$\dim\gg_{\da+\lam}=1$, then $\a:=\da+\lam\in\rtimes.$ So if
$e_\a$, $f_\a$, $h_\a$ are as in (\ref{sl2}), then for any
$x\in\gg_{\db+\mu}$ $(\da\in\rd, \mu\in\Lam)$
$$
[h_\a,x]=(\db+\mu)(h_\a)x=\big(t_{\db+\mu},\frac{2t_\a}{(\a,\a)}\big)x=
\frac{2(\db,\da)}{(\da,\da)}x=(\db,\da^\vee)x.
$$
Thus part (ii) of LT2 also holds for $\gg$.

Recall that the core $\gc$ of $\gg$ is the ideal of $\gg$
generated by root spaces $\gg_\a$, $\a\in\rtimes= R_t^\times.$  So
$\gc$ inherits from $\gg$ a $(\rd\oplus\Lam)$-grading, namely
$$
\gc=\bigoplus_{\a\in\lrd\oplus\Lam}(\gc)_\a,\quad\hbox{where}\quad
(\gc)_\a=\gc\cap\gg_\a.
$$
Moreover, from the way $\gc$ is defined, we have
\begin{equation}\label{gcd}
(\gc)_\d=\sum_{\a\in\rtimes}[\gg_{\a+\d},\gg_{-\a}]=\sum_{\da\in\rd^\times}
\sum_{\lam\in{\Lam_t}}[\gg_{\da+\lam+\d},\gg_{-\da-\lam}]\quad\quad(\d\in
R^0).
\end{equation}

Next let $\gt=\gc/{\mathcal Z}(\gc)$.  Then
$$
\gt=\bigoplus_{\a\in\lrd\oplus\Lam}\gt_\a\quad\hbox{where}\quad\gt_\a=\frac{(\gc)_\a+{\mathcal
Z}(\gc)}{{\mathcal Z}(\gc)},
$$
is a $\lrd\oplus\Lam$-grading for $\gt$. Note that if $\a\in\riso$
then by Proposition \ref{w6}(ii), $\gt_\a\sub\zgc$ and so
$\gt_\a=\{0\}$. Therefore we may assume that $\a\in R_t$. Thus
\begin{equation}\label{gt}
\quad\gt=\bigoplus_{\a\in\lrd\oplus\Lam_t}\gt_\a.
\end{equation}
That is $\gt$ has a $\lrd\oplus\Lam_t$-grading. Clearly, we have
$\gt_{\da+\lam}=\{0\}$ if $\da\not\in\rd.$  Thus LT1 holds for
$\gt$ with respect to the $(\lrd\oplus\Lam_t)$-grading (consider
$\ll_{\da}^{\lam}=\gt_{\da+\lam},\da\in\lrd,\lam\in\Lam_t$).

For $\da\in\rd^\times$ and $\lam\in\Lam_t$,
$\dim\gt_{\da+\lam}\leq 1$, as $\dim\gg_{\da+\lam}\leq 1$.
Moreover, by Theorem \ref{aabgp}(d) and Proposition \ref{w6},
$\gg_{\da}\cap\zgc=\{0\}$ and so
\begin{equation}\label{dim}
\dim\gt_{\da}=1\quad\hbox{for}\quad\da\in\rd_{\ind}\sub \rtimes.
\end{equation}
Thus part (i) of LT2 holds for $\gt$. Now as part (ii) of LT2
holds for $\gg$, one can see that it also holds for $\gt$ by
considering
$$
\tilde{e}_{\a}=e_\a+\zgc,\quad\tilde{f}_\a=f_\a+\zgc\andd\tilde{h}_\a=h_\a+\zgc.
$$
From (\ref{gcd}) we see that LT3 holds for $\gt$.

Finally, we show that LT4 holds for $\gt$. So let $\d\in R^0_t$.
Then $\gg_\d\not=\{0\}$ and $\d$ is not isolated. Thus by
Proposition \ref{w6}, $\gg_\d\not\sub\zgc$. So we have
$\gt_\d\not=\{0\}$ and
$$
\d\in\la\d\in\Lam_t\mid\gt_{\da+\lam}\not=\{0\}\hbox{ for some
}\d\in\rd\ra.
$$
This shows that LT4 holds for $\gt$ and completes the proof.\qed

\begin{cor}\label{final}
The core modulo center of a GRLA is a direct sum of centerless Lie
tori.
\end{cor}

We remark here that if $\gg$ is an EALA, then Theorem
\ref{classification} is a consequence of [AG, Proposition 1.28].
Also an statement similar to Theorem \ref{classification} is
announced in [N2, Proposition 3] for a class of Lie algebras which
includes the class of EALAs, however GRLAs do not necessarily
satisfy GR6(b), while the Lie algebras appearing in [N2] are tame
by definition and so satisfy GR6b. In [N2] a procedure is
introduced for the construction of an EALA starting from a
centerless Lie torus. In fact it is shown that all EALAs arise
this way. It is therefore natural to ask if one can introduce a
similar procedure for constructing a GRLA starting from a direct
sum of centerless Lie tori.

\vspace{5mm}
\section{\bf Construction of new GRLAs from old}
\setcounter{equation}{0} It is known that affine Lie algebras can
be realized by a process known as affinization-and-twisting [K].
It is also known that affine Lie algebras can be realized as the
fixed points of some others under a finite order automorphism.
This phenomenon has recently been investigated for the class of
EALAs (see [ABP], [ABY] and [A2]). In this section we consider a
similar method for constructing new GRLAs from old ones.

Let $\gfh$ be a GRLA with root system $R$. Let $\sg$ be an
automorphism of $\gg$ and set
$$
\gs=\{x\in\gg\mid\sg(x)=x\}\andd\hs=\{h\in\hh\mid\sg(h)=h\},
$$
that is $\gs$ and $\hs$ are fixed points of $\gg$ and $\hh$ under
$\sg$, respectively. Assume that $\gg$ satisfies
\begin{equation}\label{auto}
\begin{array}{l}
\bullet\;\; \sg \hbox{ is of finite order},\\

\bullet\;\; (\sg x,\sg y)=(x,y) \hbox{ for all }x,y\in\gg,\\

\bullet\;\; \sg(\hh)=\hh,\\

\bullet\;\; \hbox{ The centralizer of }\hs \hbox{ in $\gs$ equals
$\hs$.}
\end{array}
\end{equation}

\begin{thm}\label{thm1}
Let $\gfh$ be a GRLA and $\sg$ be an automorphism of $\gg$
satisfying (\ref{auto}). Then $(\gs,\fm,\hs)$ satisfies GR1-GR4,
where $\fm$ is the form on $\gg$ restricted to $\gs$. In
particular, $\gs$ is a GRLA if its root system contains some
nonisotropic roots.
\end{thm}

\proof It is shown in [ABY] that if $\gfh$ is an EALA then $\gs$
satisfies GR1-GR4. Now checking the proof of [ABY, Theorem 2.63],
one can see that the irreducibility of $\gg$ (or its root system
$R$) is not used at all to prove that $\gs$ satisfies
GR1-GR4.\hfill\qed\vspace{3mm}

 Next we consider the
so called {\it affinization} of  a Lie algebra $\gg$ introduced in
[ABP]. Let $\gg$ be a complex Lie algebra and let $c$ and $d$ be
two symbols. Consider the vector space
$$\aff(\gg):=(\gg\otimes{\bbbc}[t,t^{-1}])\oplus{\bbbc}c\oplus{\bbbc}d,
$$
where $\bbbc[t,t^{-1}]$ is the algebra of Laurent polynomials in
variable $t$. Then $\aff(\gg)$ becomes a Lie algebra by assuming
that $c$ is central, $d=t\frac{d}{dt}$ is the degree derivation so
that $[d,x\otimes t^n]=nx\otimes t^n$, and
$$ [x\otimes t^n,y\otimes t^m]=[x,y]\otimes
t^{n+m}+n(x,y)\d_{m+n,0}c.$$ $\aff(\gg)$ is called the {\it
affinization} of $\gg$.

If $\gg$ is equipped with a invariant symmetric bilinear form,
then one can define an invariant symmetric bilinear form
$\fm_{\aff}$ on $\aff(\gg)$ by
%\begin{equation}%\label{form1}
$$
(\a x\otimes t^n+\b c+\gamma d,\a'y\otimes t^m+\b'c+\gamma'
d))_{\aff}=\a\a'\d_{n,-m}(x,y)+\b\gamma'+\b'\gamma.
$$
%\end{equation}
Moreover, this form is nondegenerate if the form on $\gg$ is
nondegenerate.
\begin{thm}\label{thm2}
Let $\gfh$ be a GRLA with root system $R$ and let $\sg$ be an
automorphism of $\gg$ satisfying (\ref{auto}). Then
$$\big(\aff(\gg),\fm_{\aff},\hh\oplus{\bbbc}c\oplus{\bbbc}d\big)$$ is a
GRLA with root system $\rt=R+{\bbbz}\d$ where $\d\in\hht^{\star}$
is defined by $\d(d)=1$ and $\d(\hh+{\bbbc}c)=0$. Moreover,
$\aff(\gg)$ is tame if and only if $\gg$ is tame. Finally if we
extend $\sg$ to an automorphism of $\aff(\gg)$ by
$$\sg(x\otimes t^i+rc+sd)=\zeta^{-i}\sg(x)\otimes t^i+rc+sd,
$$
then $\sg$ satisfies (\ref{auto}) and $\aff(\gg)^\sg$ satisfies
GR1-GR4. In particular, $\aff(\gg)^\sg$ is a GRLA if its root
system contains some nonisotropic roots.
\end{thm}

\proof It can be checked easily that $\aff(\gg)$ is a GRLA with
root system $\rt$ as in the statement, and that $\sg$ extended to
$\aff(\gg)$ satisfies (\ref{auto}). The statement regarding
tameness is also easy to see. The last statement now follows from
Theorem \ref{thm1}.\hfill\qed\vspace{3mm}

\begin{exa}\label{exam3} {\rm Let $\gg=\oplus_{i=1}^{k}\dot{\gg}_i$ be a
complex semisimple Lie algebra with $k>1$, where each
$\dot{\gg}_i$ is a simple Lie algebra and consider $\aff(\gg)$.
According to Theorem \ref{thm2}, $\aff(\gg)$ is a GRLA. Now if we
follow the same procedure as in the proof of Theorem \ref{main},
we see that
$$\aff(\gg)=\sum_{i=1}^{k}\gg_i,\quad\hbox{where}\quad\gg_i=\dot{\gg}_i\otimes
\bbbc[t,t^{-1}] \oplus\bbbc c\oplus\bbbc d.$$ In particular, the
sum is not direct. Note that each $\gg_i$ is an affine Lie
algebra.}
\end{exa}

\end{document}